\documentclass[12pt,oneside]{amsart}

\usepackage{amscd,amsmath,amsfonts,amssymb,amsthm}

\usepackage[all]{xy}

\newcommand{\cone}[1]{\ensuremath{\mathcal{C}(#1)}}

\newcommand{\rid}{\ensuremath{/\!\!/}}
\newcommand{\hrid}{\ensuremath{/\!\!/\!\!/}}
\newcommand{\jrid}{\ensuremath{\hrid_{\text{Joyce}}}}

\newcommand{\C}{\ensuremath{\mathbb C}}
\newcommand{\Ha}{\ensuremath{\mathbb H}}

\newcommand{\R}{\ensuremath{\mathbb R}}

\newcommand{\Z}{\ensuremath{\mathbb Z}}

\newcommand{\CP}[1]{\ensuremath{{\C\mathbb P}^{#1}}}
\newcommand{\aut}[1]{\ensuremath{\text{\upshape \rmfamily Lck}(#1)}}
\newcommand{\sas}[1]{\ensuremath{\text{\upshape \rmfamily Sas}(#1)}}

\newcommand{\lchk}[1]{\ensuremath{\text{\upshape \rmfamily Lchk}(#1)}}

\newcommand{\Hom}[1]{\ensuremath{\text{\upshape \rmfamily Hom}(#1)}}

\newcommand{\Isom}[1]{\ensuremath{\text{\upshape \rmfamily Isom}(#1)}}

\newcommand{\un}[1]{\ensuremath{\text{\upshape \rmfamily U}(#1)}}

\newcommand{\SL}[1]{\ensuremath{\text{\upshape \rmfamily SL}(#1,\Z)}}
\newcommand{\Sp}[1]{\ensuremath{\text{\upshape \rmfamily Sp}(#1)}}
\newcommand{\Tor}{\ensuremath{\text{\upshape \rmfamily Tor}}}

\newcommand{\iso}{\ensuremath{\simeq}}

\newcommand{\lck}{\ensuremath{[(K,\Gamma)]}}

\newcommand{\Kmin}{\ensuremath{K_{\text{min}}}}
\newcommand{\Kminprime}{\ensuremath{K_{\text{min}}'}}
\newcommand{\Wmin}{\ensuremath{W_{\text{min}}}}
\newcommand{\Gammamin}{\ensuremath{\Gamma_{\text{min}}}}
\newcommand{\Gammaminprime}{\ensuremath{\Gamma_{\text{min}}'}}
\newcommand{\presmin}{\ensuremath{(\Kmin,\Gammamin)}}
\newcommand{\presminprime}{\ensuremath{(\Kminprime,\Gammaminprime)}}
\newcommand{\momentmin}{\ensuremath{\moment_{\text{min}}}}

\newcommand{\Kmax}{\ensuremath{\tilde{K}}}
\newcommand{\Kmaxprime}{\ensuremath{\tilde{K}'}}
\newcommand{\Wmax}{\ensuremath{\tilde{W}}}
\newcommand{\Gammamax}{\ensuremath{\tilde{\Gamma}}}
\newcommand{\Gammamaxprime}{\ensuremath{\tilde{\Gamma}'}}
\newcommand{\presmax}{\ensuremath{(\Kmax,\Gammamax)}}
\newcommand{\presmaxprime}{\ensuremath{(\Kmaxprime,\Gammamaxprime)}}
\newcommand{\momentmax}{\ensuremath{\tilde{\moment}}}

\newcommand{\Homot}[2]{\ensuremath{\text{\upshape \rmfamily Hom}_{#1}(#2)}}
\newcommand{\Hommin}{\ensuremath{\Homot{\Gammamin}{\Kmin}}}
\newcommand{\Hommax}{\ensuremath{\Homot{\Gammamax}{\Kmax}}}

\newcommand{\lee}{\ensuremath{\omega}}

\newcommand{\fund}{\ensuremath{\Omega}}
\newcommand{\map}[3]{\mbox{${#1}\colon{#2}\to{#3}$}}
\newcommand{\liealg}[1]{\ensuremath{\mathfrak{#1}}}
\newcommand{\st}{\ensuremath{\text{ such that }}}
\newcommand{\dismap}[5]{
\[
\begin{array}{rcll}
#1: & #2 & \longrightarrow & #3\\
 & #4 & \longmapsto & #5
\end{array}
\]
}

\newcommand{\Chi}{\ensuremath{\mathfrak X}}

\newcommand{\gtwo}{\ensuremath{\text{\upshape \rmfamily G}_2}}

\newcommand{\spin}[1]{\ensuremath{\text{\upshape \rmfamily Spin}(#1)}}
\newcommand{\su}[1]{\ensuremath{\text{\upshape \rmfamily SU}(#1)}}

\renewcommand{\Im}[1]{\ensuremath{\text{\upshape \rmfamily Im}#1}}

\newcommand{\moment}{\ensuremath{\mu}}

\numberwithin{equation}{section}

\newtheorem{te}{Theorem}[section]
\newtheorem*{te*}{Theorem}
\newtheorem{pr}[te]{Proposition}
\newtheorem{co}[te]{Corollary}

\theoremstyle{definition}
\newtheorem{de}[te]{Definition}
\newtheorem{re}[te]{Remark}
\newtheorem{ex}[te]{Example}

\newcommand{\EndDim}{\ensuremath{\nopagebreak\hfill\blacksquare}}
\newenvironment{D}[1][]{{\nopagebreak\noindent\em Proof#1: }}{\EndDim}

\begin{document}

\title[Reduction of Vaisman structures]{Reduction of Vaisman structures\\ in
complex and quaternionic
geometry
}

\author{Rosa Gini, Liviu Ornea, Maurizio Parton, Paolo
Piccinni}
\address{Rosa Gini,\newline via Fucini 48, 56127 Pisa, Italy.}
\email{rosa.gini@poste.it}
\address{Liviu Ornea,\newline University of Bucharest, Faculty of Mathematics,
14 Academiei str.,
70109 Bucharest, Romania.}
\email{Liviu.Ornea@imar.ro}
\address{Maurizio Parton,\newline Universit\`a di Chieti--Pescara\\
Dipartimento di Scienze, viale Pindaro $87$, I-65127 Pescara, Italy.}
\email{parton@sci.unich.it}
\address{Paolo Piccinni,\newline Universit\`a degli Studi di Roma "La Sapienza"\\
Dipartimento di Matematica "Guido Castelnuovo", Piazzale Aldo Moro 2, I-00185,
Roma, Italy.}
\email{piccinni@mat.uniroma1.it}

\begin{abstract}
%
We consider locally conformal K\"ahler geometry as an equivariant (homothetic) K\"ahler geometry: a locally conformal K\"ahler manifold 
is, up to equivalence, a pair
$(K,\Gamma)$ where $K$\/ is a K\"ahler manifold and $\Gamma$ a discrete Lie group of biholomorphic homotheties acting freely and properly discontinuously. We define  a new invariant of a locally conformal K\"ahler manifold $(K,\Gamma)$ as the rank of a natural quotient of $\Gamma$, and prove its invariance under reduction.
This equivariant point of view leads to a proof that locally conformal K\"ahler reduction of 
compact Vaisman manifolds produces Vaisman 
manifolds and is equivalent to a Sasakian reduction. 
Moreover we define locally conformal hyperk\"ahler reduction as an equivariant version of hyperk\"ahler reduction and
in the compact case we show its equivalence with 3-Sasakian reduction. Finally we show that locally 
conformal hyperk\"ahler reduction induces hyperk\"ahler
with torsion (HKT) 
reduction of the associated HKT structure and the 
two reductions are compatible, even though not 
every HKT reduction comes from a locally conformal 
hyperk\"ahler reduction.
\end{abstract}

\maketitle



\tableofcontents

\section{Introduction}

Symplectic reduction was already extended to many other structures defined 
by a closed form. Among them, recently arrived in the field, is the reduction of locally conformal
K\"ahler structures, see~\cite{GOPLCK}. 
One of the main results in the above paper concerned the conditions under 
which a particular class of compact locally conformal K\"ahler manifolds, 
namely Vaisman manifolds, is preserved  by reduction.  

Compact Vaisman manifolds are, in a certain way, equivalent to Sasakian 
manifolds, as their universal cover is a Riemannian cone over a Sasakian manifold. 
This was implicit in \cite{KaOGFC} and was made explicit in \cite{GOPLCK}. 
On the other hand, the Structure Theorem in \cite{OrVSTC} proves that any 
compact Vaisman manifold is a Riemannian suspension over a circle, with 
fibre a compact Sasakian manifold. 

It is one of the purposes of the present paper to make clear the distinction 
between these two results in the language of \emph{presentations}.
More precisely, in \cite{GOPLCK} the authors verified that considering a locally conformal K\"ahler 
manifold as a pair $(K,\Gamma)$, $K$\/ K\"ahler and $\Gamma$ acting by biholomorphic homotheties, 
provides an interesting insight in locally conformal K\"ahler 
and Vaisman geometry. Those are considered as equivariant versions of K\"ahler and Sasakian geometry, considered as \emph{homothetic}\/ geometries. The ``distance'' from the K\"ahler (Sasakian) geomety of $K$\/ is measured by a new invariant, the \emph{rank}\/ of the manifold $(K,\Gamma)$, that we define as the rank of the non-isometric part of the action of $\Gamma$.  
In this setting a Vaisman manifold is a pair $(K,\Gamma)$ where $K=\cone{W}$ is the K\"ahler cone of a Sasakian manifold and $\Gamma$ commutes elementwise with the radial flow.

In this paper we aim at delimiting how much this language can be 
considered to be exhaustive, that is, to answer the question: can locally conformal K\"ahler 
and Vaisman geometry of a manifold be completely described by a presentation of its?

The answer is generally yes, with one proviso:
 among the equivalent presentations defining a same manifold, some presentations are ``more equal'' 
than others, namely, the ``big\-gest'' and the ``smallest'' one.

The use of presentations leads to easy proofs of several interesting 
properties. In particular we 
completely drive locally conformal K\"ahler reduction, as 
defined in \cite{GOPLCK}, to K\"ahler reduction, and, in the compact Vaisman case, to Sasakian reduction, which in turn proves that compact Vaisman manifolds are closed under reduction.
Indeed, here we prove the following Theorem.

\smallskip

\noindent{\bf Theorem.}
{\em
The locally conformal K\"ahler reduction of a compact Vaisman manifold
is a Vaisman manifold.
}
\smallskip

This result depends heavily on a rather striking 
and more general one, namely: all the isometries of the cone metric over a compact 
Riemannian manifold project to identities of the generator lines. Here is where the compacity of the Vaisman manifold is used: one needs it to know that the minimal presentation is a cone over a {\it compact} Sasakian manifold.

Vaisman manifolds appear naturally also in quaternionic geometry: every compact 
locally conformal hyperk\"ahler manifold is Vaisman, see \cite{PPSEWE}. 
Beside their interest as a living structure in Hermitian quaternionic geometry 
(see, for instance, the classification in \cite{CaSAHS}), the interest for 
locally conformal hyperk\"ahler manifolds is also motivated by their close 
relation to 3-Sasakian structures, see \cite{OrPLCK}, and to HKT structures, 
see \cite{OPSPFH}. 
It is then only natural to place Vaisman reduction in the 
context of locally conformal hyperk\"ahler reduction, which we define and study 
 in the language of presentations, 
this  yielding to an easy construction of locally conformal hyperk\"ahler reduction: 
this is completely induced by hyperk\"ahler reduction, hence the interaction with 
Joyce's hypercomplex reduction and with HKT reduction is naturally managed.

\smallskip

The structure of the paper is as follows: each of the three main contexts (locally conformal K\"ahler geometry, Vaisman geometry, locally conformal hyperk\"ahler geometry) is analysed from two points of view: first the definition and properties of presentations, then the definition and properties of maps, hence of reduction: the approach is categorical, even though the language of category is never explicitly mentioned. Hence in Section~\ref{vaisman} the notion of presentation of a locally conformal K\"ahler manifold is defined, and accordingly a new equivalent definition of locally conformal K\"ahler manifold is introduced, together with the notion of 
{\em rank}. In Section~\ref{lckreduction} locally conformal K\"ahler automorphisms and twisted Hamiltonian actions are restated in terms of presentations, and reduction is reconstructed as an equivariant process, with some detail, also leading to the fact that rank is closed under reduction. In Section~\ref{vaipresentation} presentations of Vaisman manifolds are explained, and it is proven that compact Vaisman manifolds all have rank 1. In Section~\ref{vaireduction} Vaisman automorphisms (as defined in \cite{GOPLCK}) are characterized in the compact case, thus leading to the main Theorem that the reduction of compact Vaisman manifolds is always Vaisman: in order to obtain this statement several interesting results are provided, linking even more clearly compact Vaisman geometry with Sasakian geometry ({\bf Proposition~\ref{claimnoncompact}, Theorems~\ref{claimcompact}, \ref{maintheorem}, \ref{lastth}}). 
In Section~\ref{lchkpresentation} the language is shifted to the hypercomplex context, where some more rigidity appears, leading to a sufficient condition for a hyperk\"ahler manifold to be a cone on a 3-Sasakian manifold ({\bf Theorem~\ref{hypcone}});
in Section~\ref{lchkred} locally conformal hyperk\"ahler reduction is defined, in pure equivariant terms ({\bf Theorem~\ref{teonostro}}); finally in Section~\ref{hktred} the link with Joyce reduction and HKT geometry are explored.

\smallskip

\noindent{\bf Acknowledgements.} R.G.~thanks Paul Gauduchon for his warm interest on the language of presentations.
L.O.\ is grateful to Alex Buium, James Montaldi and Juan-Pablo Ortega for very useful discussions about group actions on manifolds. All the authors are grateful to Charles Boyer and Kris Galicki for useful remarks on Theorem~\ref{hypcone}, to Andrew Swann and to Misha Verbitsky for very pertinent observations on the first version of this paper.
 
\section{Presentations of 
locally conformal K\"ahler  manifolds}\label{vaisman}

Usually, a locally conformal K\"ahler manifold is a conformal
Hermitian manifold $(M,[g])$ such that $g$ is conformal to local K\"ahler
metrics. The conformal class $[g]$ corresponds to a unique 
cohomology class $[\lee_g]\in H^1(M)$, whose representative 
$\lee_g$ is defined as the unique closed 1-form satisfying 
$d\fund_g=\lee_g\wedge\fund_g$, where $\fund_g$ denotes the fundamental form 
of $g$. The 1-form $\lee_g$ is called the Lee form of $g$.

Let $M$ be a complex manifold of complex dimension at
least two.  If there exists a complex covering space $\tilde{M}$
of $M$ admitting a
K\"ahler metric ${g}_K$ and such that $\pi_1(M)$
acts on it by holomorphic homotheties, then one obtains a 
locally conformal K\"ahler structure $[g]$ on $M$ by pushing forward locally ${g}_K$
and glueing the local metrics via a partition of unity. More explicitly, in the conformal class $[{g}_K]$
there exists a $\pi_1(M)$-invariant metric $e^{f}g_K$, which induces a locally conformal K\"ahler metric 
on $M$. 

Conversely, let $(M,[g])$ be a locally conformal K\"ahler manifold.
The pull-back of any Lee form to the universal covering $\tilde{M}$
is exact, say $\tilde{\lee}_g=df$, and $e^{-f}\tilde{g}$ turns out to be a K\"ahler 
metric on $\tilde{M}$, such that $\pi_1(M)$ acts on it by holomorphic homotheties.
According to the fact that the K\"ahler metric $e^{-f}\tilde{g}$ is defined up to 
homotheties, this paper is concerned with homothety classes of K\"ahler manifolds. In a homothetic K\"ahler manifold $K$ we denote $\Hom{K}$ the group of biholomorphic homotheties, and
\[
\map{\rho_K}{\Hom{K}}{\R^+}
\]
the group homomorphism associating to a homothety its dilation factor, that is, defined by $f^*g=\rho_K(f)g$.
Remark that both $\rho_K$ and $\Isom{K}=\ker{\rho_K}$ are well-defined for a homothetic K\"ahler manifold $K$.

\begin{re} 
Note that this puts rather severe restrictions on the groups that can be fundamental groups of locally conformal K\"ahler manifolds: there are groups, for instance the ones generated by elements of finite order, as $\SL{2}$, which do not send any non-trivial homomorphism to $\R^+$.
\end{re}

What above means that every locally conformal K\"ahler
manifold can be obtained by a homothetic 
K\"ahler manifold $K$ and a 
discrete Lie group $\Gamma\subset\Hom{K}$ acting freely and properly
discontinously on $K$. This motivates 
the following Definition.

\begin{de}
A pair $(K,\Gamma)$ is a {\em presentation} 
if $K$ is a homothetic K\"ahler manifold and $\Gamma$ a 
discrete Lie group of biholomorphic homotheties acting freely and properly
discontinously on $K$. If $M$ is a locally conformal K\"ahler manifold and $M=K/\Gamma$ 
(as locally conformal K\"ahler manifolds, that is, they are conformally and biholomorphically
equivalent), we say that $(K,\Gamma)$ is a presentation {\em of} $M$.
\end{de}

\begin{ex}\label{hopfsurfaces}
The most classical examples of locally conformal K\"ahler manifolds
are the Hopf surfaces. Following Kodaira, a Hopf surface is any compact complex
surface with universal
covering $\C^2\setminus 0$.
Since
any discrete group of
biholomorphisms of $\C^2\setminus 0$ is of the form $H\rtimes\Z$, where $H$ is a
finite group (see \cite{KodSC1,KodSC2,KodSC3,KodSC4}), any locally conformal K\"ahler metric on a Hopf surface
can be given as a K\"ahler metric on $\C^2\setminus 0$ such that $H\rtimes\Z$ acts by
homotheties. In this setting, $H\rtimes\Z$ acts by
homotheties for $dz_1\otimes
d\bar{z}_1+dz_2\otimes d\bar{z}_2$ if and only if $H$ is trivial (that is, if the Hopf 
surface is primary), the action of $\Z$ is diagonal (that is, if the Hopf 
surface has K\"ahler rank $1$), and the diagonal action $(z_1,z_2)\mapsto(\alpha
z_1,\beta z_2)$ is given by complex numbers $\alpha$, $\beta$ of same lenght
$|\alpha|=|\beta|$. All other cases has been settled in \cite{GaOLCK} and
\cite{BelMSN}.
\end{ex}

\begin{ex}\label{newexample}
A very recent example of locally conformal K\"ahler manifold is described in \cite{OeTNKC}. For any integer
$s\ge 1$, the authors construct 
a group $\Gamma$ acting freely and properly discontinously by biholomorphisms on $H^s\times \C$ (where $H$ denotes the upper-half complex 
plane). The K\"ahler form such that $\Gamma$ acts by homotheties is $i\partial\bar{\partial}F$, where the K\"ahler potential $F$ on $H^s\times \C$  is given by
\[
F(z_1,\dots,z_{s+1})=\frac{1}{\prod_{j=1}^si(z_j-\bar{z}_j)}+|z_{s+1}|^2.
\]
Note that for $s=1$ this construction recovers the Inoue surfaces $S_M$ described in \cite{InoSCV}.
\end{ex}

\begin{de}
Let $(K,\Gamma)$ be a presentation.
The associated {\em maximal} presentation is \presmax,
where $\Kmax$ is the homothetic K\"ahler universal covering of $K$ and 
$\Gammamax$ is the lifting of $\Gamma$ to $\Kmax$.
Then 
\[
\presmin=\left(\frac{\Kmax}{\Isom{\Kmax}\cap\Gammamax},\frac{\Gammamax}{\Isom{\Kmax}\cap\Gammamax}\right)
\]
is called the associated {\em minimal} 
presentation.
\end{de}

\begin{re}
The maximal and minimal presentations associated to $(K,\Gamma)$ depend only on 
the locally conformal K\"ahler manifold $K/\Gamma$. Moreover 
$\Gammamax=\pi_1(K/\Gamma)$ and $\Gammamin=\rho_K(\Gamma)$.
\end{re}

The presentation of a locally conformal K\"ahler manifold is of
course not unique, but the maximal and minimal presentations are unique
as $\Gammamax$, $\Gammamin$-spaces respectively. This means that they can
be used to distinguish between different locally conformal K\"ahler manifolds,
as suggested by the following two Definitions.

\begin{de}
Let $(K,\Gamma)$, $(K',\Gamma')$ be two presentations. We say that 
$(K,\Gamma)$ is equivalent to $(K',\Gamma')$
if $\presmax=\presmaxprime$. The $=$ sign means that there exists an equivariant map,
namely, there exists
a pair of maps $(f,h)$,  where \map{f}{\Kmax}{\Kmaxprime} is a biholomorphic homothety, 
\map{h}{\Gammamax}{\Gammamaxprime} is an isomorphism and for any $\gamma\in\Gammamax$,
$x\in\Kmax$ we have
$f(\gamma x)=h(\gamma)f(x)$.  
Equivalently, we say that 
$(K,\Gamma)$ is equivalent to $(K',\Gamma')$
if  
$\presmin=\presminprime$.
\end{de}

\begin{de}
A {\em locally conformal K\"ahler} manifold is an equivalence class $\lck$ of presentations. 
\end{de}

It is worth to stress out once again that this latter Definition is equivalent with the usual one.
Hence, whenever we say ``[\dots]$M$ is a locally conformal K\"ahler manifold[\dots]'' the 
reader can indifferently think at the usual definition or at the above new one. Nevertheless,
one should remember that in the spirit of what is written in the Introduction, 
in this paper we always (henceforth) have the new definition in mind. 

\begin{re}\label{gck}
A locally conformal K\"ahler manifold $\lck$ is globally conformal K\"ahler if and only 
if $\Gamma\subset\Isom{K}$.
\end{re}

For any presentation $(K,\Gamma)$, since $\rho_K(\Gamma)$ is a finitely 
generated
subgroup of $\R^+$, it is isomorphic to $\Z^r$, for a certain $r$ a priori depending on the presentation $(K,\Gamma)$. The next 
Proposition shows that $r$ is actually an invariant of the locally conformal K\"ahler structure. This invariant encodes 
the ``locally conformal'' part of the structure, the K\"ahler part being encoded by $\ker\rho_K$. 

\begin{pr}
For any presentation $(K,\Gamma)$, the rank of the free abelian group 
$\rho_\Gamma(\Gamma)$  
depends only on the equivalence class $\lck$.
\end{pr}

\begin{D}
%
The maximal presentation $\presmax$ depends only on $\lck$. 
Remark that $\Gamma=\Gammamax/\pi_1(K)$ and that 
$\pi_1(K)\subset\ker{\rho_{\Kmax}}$ (see Remark~\ref{gck}). 
Moreover the following diagram (where the vertical arrow is a surjection) commutes
\[
\xymatrix{
\Gammamax\ar[dr]^{\rho_{\Kmax}}\ar@{>>}[d]\\
\Gamma\ar[r]_{\rho_{K}}&\R^+\\
}
\]
and this implies that $\Im\rho_K=\Im\rho_{\Kmax}$, hence the claim.
\end{D}

\begin{de}
The {\em rank}
of a locally conformal K\"ahler manifold $M$ is the 
non-negative integer $r_M$ given by the previous Proposition.
\end{de}

\begin{pr}\label{rM} 
Let $M$ be a locally conformal K\"ahler manifold and $r_M$ its rank. Then 
\[
0\leq r_M\leq b_1(M)
\]
and $r_M=0$ if and only if $M$ is
globally conformal K\"ahler.
In particular, if $b_1(M)=1$ then $r_M=1$.
\end{pr}

\begin{D}
Since $\rho_{\Kmax}(\pi_1(M))$ is abelian, the commutator $[\pi_1(M),\pi_1(M)]$ is a subgroup
of $\ker\rho_{\Kmax}$.
Thus,
\[
\begin{split}
r_M=\text{rank}\frac{\pi_1(M)}{\ker(\rho_{\Kmax})\cap\pi_1(M)}&\leq\text{rank}\frac{\pi_1(M)}{[\pi_1(M),\pi_1(M)]}\\
&=\text{rank}(H_1(M))= b_1(M)
\end{split}
\]
Moreover, by Remark~\ref{gck}, $r_M=0$ if and only if  $M$ is
globally conformal K\"ahler.
Finally, if $b_1(M)=1$, the manifold cannot admit any globally 
conformal K\"ahler structure, hence $r_M=1$.
\end{D}

\begin{ex}
The examples $M_s=(H^s\times \C)/\Gamma$ given in \ref{newexample} and described
in \cite{OeTNKC} satisfy $\text{rank}(M_s)=s$, for any integer
$s\ge 1$. 
\end{ex}

\begin{re}
It is worth noting that $r_M=1$ does not imply $b_1=1$. We shall see in Corollary \ref{vai_rank} that all compact Vaisman manifolds have rank $1$, but one can obtain examples with arbitrarily large $b_1$ by taking induced Hopf bundles over curves of large genus in $\CP{2}$.
\end{re}

This picture of Definitions and Remarks is summarized in the 
following diagram, which also clarifies the terminology maximal and minimal. The vertical arrows are covering maps via the action of isometries, the transversal arrows are covering maps via the action of homotheties, and $r=r_M$ is the rank of the locally conformal K\"ahler manifold.
\begin{equation}\label{presdiagram}
\xymatrix{
\Kmax\ar[rd]^{\Gammamax=\pi_1(M)}\ar[d]_{\pi_1(K)} & \\
K\ar[r]^{\Gamma}\ar[d]_{\Isom{K}\cap\Gamma} & M \\
\Kmin\ar[ru]_{\Gammamin=\rho_K(\Gamma)=\Z^{r}}& \\
}
\end{equation}

\section{Locally conformal K\"ahler automorphisms and reduction}\label{lckreduction}

A locally conformal K\"ahler automorphism is usually a biholomorphic
conformal map.

\begin{de}
Let $M$ be a locally conformal K\"ahler manifold.
Denote by $\aut{M}$ the Lie group of biholomorphic conformal
transformations of $M$.
\end{de}

Coherently with our setting, we would like to describe $\aut{M}$ in terms 
of presentations, and a natural candidate is given in the following Definition.

\begin{de}
Let $(K,\Gamma)$ be a presentation. 
Denote by \Homot{\Gamma}{K} the  
Lie group of $\Gamma$-equivariant biholomorphic
homotheties of $K$.
\end{de}

The problem is that the existence of a lifting to $K$ is generally not granted, apart
from the trivial case of the maximal presentation $\presmax$. It is exactly in this context
that the notion of minimal presentation shows its relevance.

\begin{pr}\label{rosa}
Let $(K,\Gamma)$ be a presentation, and
let $f$ be any homothety of $K$ commuting with $\Gamma$. Then $f$ commutes
with $\Isom{K}\cap\Gamma$.
\end{pr}

\begin{D}
Let $\gamma\in\Isom{K}\cap\Gamma$. Since $f$ commutes with $\Gamma$, it exists
$\gamma'\in\Gamma$ such that $f\circ\gamma=\gamma'\circ f$. Apply $\rho_\Gamma$ to
both members to obtain that $\gamma'$ lies actually in $\Isom{K}\cap\Gamma$.
\end{D}

\begin{co}\label{rosacor}
Let $M$ be a locally conformal K\"ahler manifold and $\presmin$ its minimal presentation.
Any map $f\in\aut{M}$ lifts to a map $f_\text{min}\in\Homot{\Gammamin}{\Kmin}$.
\end{co}

\begin{D}
Any biholomorphic conformal
transformation $f$ of $M$ lifts to a biholomorphic homothety $\tilde{f}$ of $\Kmax$.
Then by Proposition~\ref{rosa} the map $\tilde{f}$ commutes with $\Isom{\Kmax}\cap\Gammamax$,
and induces a biholomorphic homothety $f_\text{min}$ of $\Kmin=\Kmax/\Isom{\Kmax}\cap\Gammamax$
which is $\Gammamin$-equivariant, since it is a lifting of $f$.
\end{D}

\begin{co}\label{corpi1}
Let $M$ be a locally conformal K\"ahler manifold. Then any map 
$f\in\aut{M}$ preserves the subgroup of $\pi_1(M)$ acting on $\Kmax$ by isometries.
\end{co}

 \begin{D}
Denote by $p$ the covering associated to the minimal presentation. In view of Corollary~\ref{rosacor},
$f$ lifts to the minimal presentation, hence 
\[
f_*(p_*(\pi_1(\Kmin)))=p_*(\pi_1(\Kmin)).
\]
On the other hand $p_*(\pi_1(\Kmin))$ is the group whose action on \Kmax\ produces \Kmin, and 
by definition this is the subgroup of $\Gammamax=\pi_1(M)$ acting on $\Kmax$ by isometries.
 \end{D}

The following diagram pictures this setting and adopts a convention we are trying to stick
to, that is, we use $\tilde{\square}$ when we refer to something related 
to the maximal presentation, and $\square_\text{min}$ for something related to the
minimal presentation. 

\begin{equation}\label{lckdiagram}
\xymatrix{
\Kmax\ar[rrr]^{\tilde{f}}\ar[rd]^{\tilde{\pi}_\text{min}}\ar[rdd]_{\tilde{\pi}} & & & \Kmax\ar[dl]_{\tilde{\pi}_\text{min}}\ar[ddl]^{\tilde{\pi}} \\
& \Kmin\ar[r]^{f_\text{min}}\ar[d]^{\pi_\text{min}} & \Kmin\ar[d]_{\pi_\text{min}} & \\
& M\ar[r]_{f} & M & \\
}
\end{equation}
In other words, 
\[
\aut{M}=\Hommax/\Gammamax=\Hommin/\Gammamin.
\]
Moreover,
any Lie subgroup $G$ of $\aut{M}$
yields Lie subgroups $\tilde{G}$, $G_{\text{min}}$ of \Hommax, \Hommin\ 
respectively, such that $G=\tilde{G}/\Gammamax=G_{\text{min}}/\Gammamin$.

 
\begin{de}
A connected
Lie subgroup $G$ of \aut{M} is {\em twisted Hamiltonian} if the identity component $\tilde{G}^\circ$ of
$\tilde{G}$ is Hamiltonian on \Kmax. Equivalently, $G$ is  twisted Hamiltonian if the identity 
component $G_\text{min}^\circ$ of
$G_\text{min}$ is Hamiltonian on \Kmin.
Accordingly, we say
that a map \map{\moment}{\liealg{g}}{C^\infty(M)} (or equivalently a map
\map{\moment}{M}{\liealg{g}^*}) is a {\em momentum map} for the 
action of $G$ on $M$ if it is the quotient of the corresponding momentum map $\momentmax$ on $\Kmax$, or equivalently
if it is the quotient of the corresponding
momentum map $\momentmin$ on $\Kmin$.
\end{de}

\begin{re}
The connected components are necessary to insure the existence of the momentum maps
$\momentmax$ and $\momentmin$.
\end{re}

\begin{re}\label{unknown}
If $G$ acts freely and properly on $\moment^{-1}(0)$ then
$\tilde{G}^\circ$ acts freely and properly on $\momentmax^{-1}(0)$ and
$G_\text{min}^\circ$ acts freely and properly on $\momentmin^{-1}(0)$.
\end{re}

\begin{D} 
%
Assume that the action of a group $G$ is proper on a manifold $X$ and lifts to a covering $X'$, and let $G'$ be the
identity component of the lifted action.
Let $K'\subset X'\times X'$ be compact, and let 
\[
L'=\{(g',x')\in G'\times X'\st (g'(x'),x')\in K'\}
\]
be its preimage. We want to show that $L'$ is compact.
 
Let $(g'_i,x'_i)\in L'$ be a sequence. Then $(g'_i(x'_i),x'_i)\in K'$, hence up to subsequences there exists $(y'_0,x'_0)$ limit of $(g'_i(x'_i),x'_i)$. Then $(g_i(x_i),x_i)$ converges to $(y_0,x_0)$ and since the action of $G$ on $X$ is proper then 
\[
(g_i,x_i)\rightarrow (g_0,x_0).
\]
Hence $y_0=g_0(x_0)$, namely $(g_i(x_i),x_i)\rightarrow (g_0(x_0),x_0)$. Then there exists $g'_0\in G'$ such that $y'_0=g'_0(x'_0)$. The fact that $G'$ covers $G$ implies then that
\[
(g'_i,x'_i)\rightarrow (g'_0,x'_0)
\]
hence that the lifted action is proper. The freedom of the lifted action is straightforward.
\end{D}


This notion of twisted Hamiltonian action coincides with the one given in \cite{GOPLCK}, 
and starts a reduction process. 

In view of defining a reduction process in terms of presentations (and since
of course we would like this process to be equivalent to the one given in \cite{GOPLCK}),
we give here some technical details (marked by the \verb+\footnotesize+ environment) regarding the reduction in \cite{GOPLCK}, in terms of the 
maximal presentation \presmax. 


\bigskip

\begin{footnotesize}
First, we need some additional topological conditions on the zero set of the momentum map,
that is, we suppose $\moment^{-1}(0)$ non-empty, $0$
a regular value for \moment\
and the action of $G$ free and proper on $\moment^{-1}(0)$.

Then, some notations. Call $\tilde{\pi}$ the covering map of $\Kmax$ over $M$, whose covering maps 
are given by $\Gammamax$. Points in $\Kmax$ are denoted by $\tilde{x}$, and this of course means also that
$\tilde{\pi}(\tilde{x})=x$, where $x\in M$. Elements of $\tilde{G}^\circ$ lifting $g\in G$ are denoted by $\tilde{g}$.

Then, $\tilde{G}^\circ$ acts freely and properly on the non-empty $\momentmax^{-1}(0)$, and 
the K\"ahler reduction $\Kmax\rid \tilde{G}^\circ$ is defined. We show here that $\Gammamax$ is
$\tilde{G}^\circ$-equivariant, that is, 
\begin{equation}\label{GammaisGequiv}
\tilde{\gamma}\tilde{G}^\circ=\tilde{G}^\circ\tilde{\gamma}\qquad\text{for any }\tilde{\gamma}\in\Gammamax.
\end{equation}

Indeed, take a point $\tilde{x}\in\Kmax$, and consider the intersection 
\[
\tilde{\gamma}\tilde{G}^\circ(\tilde{x})\cap\tilde{G}^\circ\tilde{\gamma}(\tilde{x}).
\]
This is clearly closed, and it can be easily shown to be also open in both
$\tilde{\gamma}\tilde{G}^\circ(\tilde{x})$ and $\tilde{G}^\circ\tilde{\gamma}(\tilde{x})$.
Thus (remember that $\tilde{G}^\circ$
is connected), we obtain Formula~\eqref{GammaisGequiv}.

The fact that $\Gammamax$ is
$\tilde{G}^\circ$-equivariant means that $\Gammamax$ acts on 
the K\"ahler reduction $\Kmax\rid \tilde{G}^\circ$ by biholomorphic homotheties. In particular,
we have a diffeomorphism 
\[
\frac{\Kmax\rid \tilde{G}^\circ}{\Gammamax}=M\rid G.
\]
We claim that $\Gammamax$ acts properly discontinously 
on $\Kmax\rid \tilde{G}^\circ$. In fact, the quotient could not be Hausdorff otherwise.

Suppose now that there exists $\gamma\in\Gammamax$ with a fixed point $\tilde{G}^\circ\tilde{x}$ in
$\Kmax\rid \tilde{G}^\circ$. Then we can find $\tilde{g}$ such that $\gamma\tilde{x}=\tilde{g}\tilde{x}$,
and this implies $g(x)=x$ on $M$. Since the action of $G$ on $M$ is free by hypothesis, we have $g=\text{Id}_M$.
This means that $\tilde{g}$ is actually in $\Gammamax$, and since it coincides with $\gamma$ on $\tilde{x}$ it must 
be exactly $\gamma$. We have proven that 
\[
\frac{\Gammamax}{\Gammamax\cap\tilde{G}^\circ}
\]
acts freely on $\Kmax\rid \tilde{G}^\circ$.
\end{footnotesize}

\bigskip

Summing up, we have proven the following Theorem.

\begin{te}\label{presred}
Let $M$ be a locally conformal K\"ahler manifold. Let $G\subset\aut{M}$ be twisted Hamiltonian, 
and suppose that $\moment^{-1}(0)$ is non-empty, $0$ is
a regular value for \moment\
and the action of $G$ is free and proper on $\moment^{-1}(0)$.
Then 
\[
(\Kmax\rid \tilde{G}^\circ,\frac{\Gammamax}{\Gammamax\cap\tilde{G}^\circ})
\]
is a presentation of the locally conformal K\"ahler reduction $M\rid G$. 
\end{te}

Theorem~\ref{presred} works the same way if we substitute the maximal presentation with the minimal presentation.
In this case, since 
\[
G_\text{min}^\circ\cap\Gammamin=1
\]
we obtain the following Corollary.

\begin{co}
The minimal presentation for a locally conformal K\"ahler reduction $M\rid G$ is given by 
\[
(\Kmin\rid G_\text{min}^\circ,\Gammamin).
\]
\end{co}

In particular, this implies that the rank is preserved under reduction.

\begin{co}
Let $M$ be a  locally conformal K\"ahler manifold, and $G$ a twisted Hamiltonian subgroup
of $\aut{M}$. Then 
\[
r_{M\rid G}=r_M.
\]
\end{co}

This discussion motivates the hypothesis in the following Definition.

\begin{de}
Let $[(K,\Gamma)]$ be a locally conformal K\"ahler manifold. Let $G$ be a connected Hamiltonian subgroup 
of $\Homot{\Gamma}{K}$ (this implies that $\Gamma$ is $G$-equivariant) with K\"ahler momentum map \map{\moment}{K}{\liealg{g}^*}.
Suppose that $\moment^{-1}(0)$ is non-empty, $0$ is
a regular value for \moment\
and the action of $G$ is free and proper on $\moment^{-1}(0)$, so that the K\"ahler reduction
$K\rid G$ is defined. Suppose moreover
that the action of $\Gamma$ on $K\rid G$ be properly discontinous, and that 
the action of $\Gamma/(\Gamma\cap G)$ on $K\rid G$ be free. Then
\[
[(K\rid G,\frac{\Gamma}{\Gamma\cap G})]
\]
is the {\em reduction} of the locally conformal K\"ahler manifold $[(K,\Gamma)]$. 
\end{de}

%
%



\section{Vaisman presentations}\label{vaipresentation}

In the category of locally conformal K\"ahler manifolds, a remarkable subclass is given by the so-called Vaisman manifolds, after the name of the author who first recognized their importance \cite{VaiLCK,VaiGHM}. They are usually defined as locally conformal K\"ahler manifolds having parallel Lee form, and this definition belongs to the Hermitian setting. 

However, it is possible to define Vaisman manifolds using the point of view of presentations. This is done in Theorem~\ref{lckcpt}, with the aid of the following Proposition
which can be proven by direct computation. 

\begin{pr}\label{claimnoncompact}
Let $W$ be a Riemannian manifold and let $f$ be a homothety
of its Riemannian cone $\cone{W}$ with dilation factor $\rho_f$. 
Suppose moreover that $f$ commutes with the radial flow $\phi_s$ given
by $\phi_s(w,t)=(w,s\cdot t)$.
Then $f(w,t)=(\psi_f(w),\rho_f\cdot t)$, where $\psi_f$ is an isometry
of $W$. In particular, the isometries of $W$ are the isometries of the cone 
$\cone{W}$ commuting with the radial flow $\phi_s$.
\end{pr}

\begin{te}\label{lckcpt}
A locally conformal K\"ahler manifold $\lck$ 
 is Vaisman if and only if $K$ is the
K\"ahler cone $\cone{W}$ of a Sasakian manifold $W$ and 
$\Gamma$ commutes elementwise with the radial flow $\phi_s(w,t)=(w,s\cdot t)$ of the cone.
\end{te}

\begin{D}
The universal covering of a Vaiman manifold is a cone over a Sasakian manifold, as was implicit in the work of Vaisman~\cite{VaiLCK} where, using the de Rham Theorem, observed the decomposition of the universal cover as the product of the real line with a Sasakian manifold that covers the leaves of the foliation ${\lee^\sharp}^\perp$, these ones having an induced Sasakian structure. That the lifted globally conformal K\"ahler metric is conformal with a cone metric is then immediate. This fact, together with Proposition~\ref{claimnoncompact}, implies that all of the presentations of a Vaisman manifold are cones.

Vice versa, it has been proven in~\cite{GOPLCK} that any manifold admitting such a presentation is Vaisman.
\end{D}

\begin{de}\label{defvai}
A pair $(\cone{W},\Gamma)$ is a {\em Vaisman presentation}  
if $W$ is a Sasakian manifold, \cone{W} its K\"ahler cone and $\Gamma$ a 
discrete Lie group of biholomorphic homotheties commuting elementwise with 
the radial flow $\phi_s$ and acting freely and properly
discontinously on $\cone{W}$. Correspondingly, the equivalence class
$[(\cone{W},\Gamma)]$ is called a Vaisman manifold.
\end{de}

Thus, whenever $M$ is a Vaisman manifold, $\Kmax$ is replaced by $\cone{\Wmax}$ and 
$\Kmin$ is replaced by $\cone{\Wmin}$, where $\Wmax$ is a simply connected Sasakian manifold and 
$\Wmin$ is the ``smallest'' Sasakian manifold such that its K\"ahler cone $\cone{\Wmin}$ covers $M$.
Moreover, the additional property of $\Gammamax$, $\Gammamin$ commuting with the radial flow
holds.
 
\begin{ex}
Consider primary Hopf surfaces of K\"ahler rank 0, that is, of the form $(\C^2\setminus 0)/\Gamma$, 
where $\Gamma$ is generated
by $(z,w)\mapsto(\alpha z +\lambda w^m,\alpha^m w)$ for a non-zero complex number $\lambda$.
We know from \cite{BelMSN} that on this kind of surfaces there are no
Vaisman structures. Accordingly, it is easy to check that $\Gamma$ does not commute with the radial 
flow of the standard cone structure on $\C^2\setminus 0$.
\end{ex}

%
%

%




\begin{ex}
The example given in \ref{hopfsurfaces} of Hopf surfaces can be generalized
to any dimension. Consider the compact complex manifold 
\[
H_{\alpha_1,\dots,\alpha_n}=\frac{\C^n\setminus0}{\Gamma}
\] 
with $\Gamma$ generated by
$(z_1,\dots,z_n)\mapsto(\alpha_1
z_1,\dots,\alpha_n z_n)$, where $\alpha_i$ are complex numbers of equal lenght $\ne 1$.
Then $\Gamma$ acts by homotheties for the standard metric $dz_1\otimes
d\bar{z}_1+\dots+dz_n\otimes d\bar{z}_n$, and the resulting $H_{\alpha_1,\dots,\alpha_n}$ is locally conformal K\"ahler.
Moreover, since $\C^n\setminus 0$ with the standard metric is the K\"ahler cone of the standard Sasakian structure 
on the sphere $S^{2n-1}$, and $\Gamma$ commutes with the radial flow,
 the Hopf manifold $H_{\alpha_1,\dots,\alpha_n}$ is Vaisman.
The case of $\alpha_i$ generic is settled in \cite{KaOGFC} deforming the standard Sasakian structure on the sphere, thus obtaining 
a Vaisman structure on any $H_{\alpha_1,\dots,\alpha_n}$.
\end{ex}
%
%
%
%

For  a compact Vaisman manifold $M$, the picture can be made more precise: one can associate a compact Sasakian manifold
$W$ as any of the fibers of a canonically defined Riemannian submersion
$M\rightarrow S^1$ and a biholomorphic homothety $\phi$ of the K\"ahler cone of
$W$ such that $M=[(\cone{W},\langle \phi \rangle)]$ 
\cite[Structure Theorem]{OrVSTC}. Note that if $(\cone{W'},\Z)$ is any other presentation
of $M$, then (up to homotheties) $W'=W$. Since any presentation of the form $(K,\Z)$ is necessarily
minimal, we have the following Corollary.

\begin{co}
The minimal presentation of any compact Vaisman manifold is $(\cone{\Wmin},\Z)$, for a compact 
Sasakian manifold $\Wmin$. 
\end{co}

\begin{co}\label{vai_rank}
If $M$ is a compact Vaisman manifold then $r_M=1$.
\end{co}

%

\section{Vaisman automorphisms and reduction}\label{vaireduction}

At this point, we are left with the description of automorphisms for Vaisman 
presentations.
Given the fact that the additional structure in the Vaisman case is the compatibility
of $\Gamma$ with the radial flow, we could define \sas{M} as those maps $f\in\aut{M}$ such that 
$\tilde{f}$ (or equivalently, as we will see, $f_\text{min}$) commutes with the radial flow:
in fact commuting with the radial flow implies reducing to the Sasakian structure of the base of the cone, due to Proposition~\ref{claimnoncompact}. 
%

Hence \sas{M} can be equivalently defined as the group of locally conformal K\"ahler automorphism induced by the Sasaki automorphisms of the base of the maximal cone via
\[
f(w,t)=(\psi_f(w),\rho_f\cdot t).
\]
The first result of this Section, in Theorem~\ref{maintheorem}, is that in the compact case the additional hypothesis on a $\Gamma$-invariant automorphism of commuting with the radial flow is actually redundant, hence $\sas{M}=\aut{M}$.

The main tool to prove the equivalence is
the following Theorem,  together with the Structure Theorem
in \cite{OrVSTC}.
\begin{te}\label{claimcompact}
For any compact
Riemannian manifold $W$ all the homotheties of the cone $\cone{W}$ are given by 
$(w,t)\mapsto(\psi(w),\rho\cdot t)$, where $\rho$ is the dilation factor and $\psi$ is an 
isometry of $W$.
In particular, all the isometries of the cone $\cone{W}$ come from
those of $W$.
\end{te}

\begin{D}
Let $d$ and
$\tilde{d}$ be the distances on $W$ and on
the cone $W\times\R^+$, respectively.
The metric completion $(W\times\R^+)^*$ of $\tilde{d}$ is
obtained by adding just one point, which we call the origin $0$. This can be
seen as follows.

Since the cone metric is given by $\tilde{g}=dr^2+r^2g$,
the length of germs of curves $[\gamma(t)]=[(w(t),r(t))]$ satisfies
\[
\tilde{g}(\dot{\gamma}(t),\dot{\gamma}(t))^{1/2}
=(\dot{r}^2(t)+r^2(t)g(\dot{w}(t),\dot{w}(t)))^{1/2}\geq |\dot{r}(t)|,
\]
and this implies the following inequality for $\tilde{d}$:
\[
\tilde{d}((w,r),(v,s))\ge |s-r|.
\]
Thus if $\{(w_n,r_n)\}$ is a Cauchy sequence not converging in
$W\times\R^+$ then 
$r_n\rightarrow 0$.
If
$(w,r)$, $(v,s)$ are points of $W\times\R^+$, then 
the 
triangular inequality
implies that
\[
\tilde{d}((w,r),(v,s))\le \min\{r,s\}d(w,v)+|s-r|.
\] 
It follows that if $\{(w_n,r_n)\}$, $\{(v_n,s_n)\}$ are Cauchy sequences 
not
converging
in
$W\times\R^+$, then $\tilde{d}((w_n,r_n), (v_n,s_n))\rightarrow 0$.

Next, observe that every
isometry of the Riemannian cone $W\times\R^+$ can be extended by mapping
$0\mapsto 0$ to a trasformation of $(W\times\R^+)^*$
preserving $\tilde{d}^*$, that is preserving rays and level submanifolds of the 
cone. This means that every isometry of the cone $W\times\R^+$  is of
the form
$(\psi,\text{Id}_{\R^+})$, where $\psi$ is an isometry of $W$. The
claim for a general homothety follows.
\end{D}

\begin{re}
Note that to prove that $r_n\rightarrow 0$ we only need the completeness of
$W$, whereas  to prove that $\tilde{d}((w_n,r_n), (v_n,s_n))\rightarrow 0$
we need its
compactness. Moreover, the compactness hypothesis is essential to complete 
the cone metric with only one point:
take for example $W=\R$, and the Cauchy sequences
$\{(-n,1/n)\}$, $\{(n,1/n)\}$. 
\end{re}
The following statement is a small {\em detour}.

\begin{co}\label{extension}
Let $\Gamma$ be a discrete Lie group of homotheties acting freely and properly
discontinously on a Riemannian cone \cone{W}, where $W$
is a compact Riemannian manifold. Then $\Gamma\iso I\rtimes\Z$, where $I$ is
a finite subgroup of isometries of $W$.
\end{co}

\begin{D}
Let \map{\rho}{\Gamma}{\R^+} be the map defined by the dilation factor of elements
of $\Gamma$.
%
Note that if $\rho(\Gamma)$ is not
cyclic, then it contains two elements $\alpha$, $\beta$ such
that $\log\alpha/\log\beta$ is irrational. This in turn implies that
$\rho(\Gamma)$ is dense in $\R^+$, a fact which, together with the compactness of $W$,
contradicts the proper discontinuity of $\Gamma$. Thus
$\rho(\Gamma)\iso\Z$. Moreover, since isometries on a Riemannian cone
$\cone{W}$ are actually isometries of the base space $W$ by Theorem~\ref{claimcompact},
the proper discontinuity of $\Gamma$
and the compactness of $W$ imply that the normal subgroup $I=\rho^{-1}(1)$ of isometries
in $\Gamma$ is finite, and $\Gamma$ is a finite extension of \Z
\[
0\rightarrow I\rightarrow\Gamma\rightarrow\Z\rightarrow 0.
\]
This implies that $\Gamma\iso I\rtimes\Z$.
\end{D}

%
%

We can now prove the first result of this section.

\begin{te}\label{maintheorem}
Let $M$ be a compact Vaisman manifold and let $f\in\aut{M}$. Let
$(\cone{W},\Gamma)$ be a Vaisman presentation of $M$ and suppose
there exists a lifting $f'\in\Homot{\Gamma}{\cone{W}}$ of $f$.
Then 
\[
f'(w,t)=(\psi_{f'}(w),\rho_{f'}\cdot t)
\]
where $\psi_{f'}$ is a Sasakian automorphism of $W$. Hence
\[
\sas{M}=\aut{M}.
\]

In particular, if $f'$ 
is an isometry then $f'=\psi_{f'}\times \text{Id}$.
\end{te}

\begin{D}
Consider the following 
diagram.

\begin{equation}
\xymatrix{
\cone{W}\ar[r]^{f'}\ar[d]_{\pi} & \cone{W}\ar[d]^{\pi} \\
\cone{\Wmin}\ar[r]^{f_\text{min}} & \cone{\Wmin}
}
\end{equation}
Since \cone{\Wmin} is K\"ahler, it follows that the covering maps of $\pi$ are isometries. 
According to Theorem~\ref{lckcpt} and Proposition~\ref{claimnoncompact}, 
$\pi$ is given by a 
projection of Sasakian manifolds. Moreover,  $f_\text{min}=\psi\times\rho$ by
Theorem~\ref{claimcompact}, and the diagram becomes the following:

\begin{equation}
\xymatrix{
\cone{W}\ar[r]^{f'}\ar[d]_{\pi\times\text{Id}} & \cone{W}\ar[d]^{\pi\times\text{Id}} \\
\cone{\Wmin}\ar[r]^{{\psi}\times\rho} & \cone{\Wmin}
}
\end{equation}
This implies 
\[
(\pi\times\text{Id})\circ f'=(\psi\times\rho)\circ(\pi\times\text{Id})
=(\psi\circ\pi)\times\rho
\]
and thus $f'(w,t)=(f'_1(w,t),\rho\cdot t)$, where $f'_1$ satisfies
\[
\pi(f'_1(w,t))=\psi(\pi(w)).
\]
This last equation implies that $f'_1(w,\R)$ is discrete in $\Wmin$, and 
\[
f'(w,t)=(f'_1(w,t),\rho\cdot t)=(f'_1(w),\rho\cdot t).
\]
To end the proof, note that $f'$ being a K\"ahler automorphism forces $f'_1$ 
to be a Sasakian automorphism.
\end{D}

The following is Diagram \eqref{lckdiagram}, adapted to the compact Vaisman case.

\begin{equation}\label{vaismandiagram}
\xymatrix{
\cone{\Wmax}\ar[rrrr]^{\psi_{\tilde{f}}\times\rho_{\tilde{f}}}\ar[rd]^{\tilde{\pi}_\text{min}\times\text{Id}}\ar[rdd]_{\tilde{\pi}} & & & & \cone{\Wmax}\ar[dl]_{\tilde{\pi}_\text{min}\times\text{Id}}\ar[ddl]^{\tilde{\pi}} \\
& \cone{\Wmin}\ar[rr]^{\psi_{f_\text{min}}\times\rho_{f_\text{min}}}\ar[d]^{\pi_\text{min}} & & \cone{\Wmin}\ar[d]_{\pi_\text{min}} & \\
& M\ar[rr]_{f} & & M & \\
}
\end{equation}

In \cite{GOPLCK} the notion of twisted Hamiltonian action for compact Vaisman manifolds was partially linked to the cone structure of the presentation. Those locally conformal K\"ahler automorphisms induced by isometries of a covering cone of the form
\[
f=\psi_f\times \text{Id}
\]
$\psi_f$ being a Sasakian isometry of the base, were called {\em Vaisman automorphisms} of the compact Vaisman manifold. It was shown that if all of the elements of a group $G\subset\aut{M}$ are Vaisman automorphisms then its action is twisted Hamiltonian.  It was moreover proven that twisted Hamiltonian actions of Vaisman automorphisms induced a Vaisman structure on the reduction. We can now prove that there are no twisted Hamiltonian actions but those ones.

\begin{te}\label{lastth}
The action of a connected group $G$ on a compact Vaisman manifold $M$ is twisted Hamiltonian if and only if it is by Vaisman automorphisms; if moreover the topological hypothesis of the Reduction Theorem are satisfied then
\[
M\rid G=[(\cone{\Wmin\rid G_\text{min}^\circ},\Gammamin)]=[(\cone{\Wmax\rid \tilde{G}^\circ},\frac{\Gammamax}{\Gammamax\cap\tilde{G}^\circ})]
\]
\end{te}

\begin{D} That the action of a group of Vaisman automorphisms admits a momentum map, 
hence is twisted Hamiltonian, is proven in Theorem 5.13 of \cite{GOPLCK}.

Vice versa let $G$ be any twisted Hamiltonian connected Lie subgroup of
$\aut{M}$. Since $\aut{M}=\sas{M}$ the elements of $\tilde{G}^\circ$ are of the form
\[
f(w,t)=(\psi_f(w),\rho_f\cdot t),
\]
and since $\tilde{G}^\circ$ is connected $\rho|\tilde{G}^\circ$ is constantly equal to 1, that is $\rho_f=1$ for all $f\in\tilde{G}^\circ$, hence $G$ is of Vaisman automorphisms.

The final claim is Theorem 5.15 of \cite{GOPLCK} in the present notation.
\end{D}

This proves the 
Theorem in the Introduction.

\begin{re}
The Theorem in the Introduction holds true for any Vaisman presentation
$[(\cone{W},\Gamma)]$ and $G$ acting on $W$. 
\end{re}
%
 
\begin{co}\label{structure}
The Structure Theorem of \cite{OrVSTC} and
Vaisman reduction are compatible.
\end{co}

\begin{D}
Note that in the equivalence 
\[
[(\cone{\Wmin},\Gammamin)]\rid G=[(\cone{\Wmin\rid G_\text{min}^\circ},\Gammamin)]
\]
the group
$\Gammamin$ is preserved. This implies that the reduced Sasakian manifold $\Wmin\rid G_\text{min}^\circ$ is
still given by the Structure Theorem. 
\end{D}

\begin{re}
Whereas the minimal presentation is preserved, as shown above in Corollary \ref{structure}, the maximal presentation is not, because $\Wmax\rid\tilde{G}^\circ$ needs not be simply connected. 
\end{re}
 
\begin{re}
If $\cone{\Wmin}\stackrel{\Z}{\rightarrow}M$ is the covering given by
the minimal presentation of a compact Vaisman manifold
and $\pi$ the projection of the cone
onto its radius, it is easy to check that $\pi$ is equivariant
with respect to the action of the covering maps on $\cone{\Wmin}$ and the
action of $n\in\Z$ on $t\in\R$ given by $n + t$.
This describes in an alternative way the projection over $S^1$ 
of the Structure Theorem in \cite{OrVSTC}, 
and moreover provides a structure theorem for other types of locally conformal structures, 
where the compactness of the base of the cone is given by other ways
(see~\cite{IPPLCP} for an application to \gtwo, \spin{7} and \spin{9} structures).
\end{re}


%




  

\section{Locally conformal hyperk\"ahler presentations}\label{lchkpresentation}

Usually, a locally conformal hyperk\"ahler manifold is a conformal
hyperhermitian manifold $(M,[g])$ such that $g$ is conformal to local
hyperk\"ahler metrics.

The same way as for locally conformal K\"ahler manifolds, we observe that a
hyperk\"ahler manifold $H$ and a discrete Lie group $\Gamma$ of hypercomplex
homotheties acting freely and properly discontinously on $H$ give rise to a 
locally conformal hyperk\"ahler manifold $H/\Gamma$, and vice versa that every 
locally conformal hyperk\"ahler manifold  can be viewed this way.

The hard part of the work being already done in the locally conformal K\"ahler setting, 
we sketch here the Definitions and properties shared in the hyperk\"ahler case. 
For the sake of simplicity, we don't use a different notation/terminology.

\begin{de}
Given a homothetic hyperk\"ahler manifold $H$ and a 
discrete Lie group $\Gamma$ of hypercomplex homotheties acting freely and properly
discontinously on $H$, the pair $(H,\Gamma)$ is called a {\em presentation}. If a locally
conformal hyperk\"ahler manifold $M$ is given, and $M=H/\Gamma$ as locally
conformal hyperk\"ahler manifolds, we say that $(H,\Gamma)$
is a presentation {\em of} $M$.
\end{de}

As before, the maximal and the minimal presentations are defined, and
gives a way to distinguish between equivalent presentations.

\begin{de}
A {\em locally conformal hyperk\"ahler} manifold is an equivalence class $[(H,\Gamma)]$ of presentations. 
\end{de}

\begin{ex}\label{quathopf}
Let $q\in\Ha\setminus 0$ be a non-zero 
quaternion of lenght different from $1$. Then the right multiplication by $q$  generates a
cocompact, free and properly discontinous action of homotheties on
$\Ha^n\setminus 0$ equipped with the standard metric. Denote by $\Gamma$ the infinite cyclic group 
$\langle h\mapsto h\cdot q\rangle$. Then $(\Ha^n\setminus 0,\Gamma)$ is a compact locally conformal 
hyperk\"ahler manifold, called a quaternionic
Hopf manifold.
\end{ex}


By definition, any locally conformal hyperk\"ahler manifold is locally
conformal K\"ahler, but much more is true in the compact case.

\begin{te}\label{hypcpt}\cite{CaPEWG,GauSWE}
A compact locally conformal
hyperk\"ahler manifold is
Vaisman for each one of its three underlying locally conformal K\"ahler
structures.
\end{te}
%

This shows that in the compact case locally conformal hyperk\"ahler
manifolds are the hypercomplex analogous of Vaisman manifolds and 
not only of locally conformal K\"ahler manifolds. This leads to the following Theorem. Before stating it, recall that a Riemannian manifold $(W,g)$ is $3$-Sasakian if, by definition, the Riemannian cone metric on $W\times \R^+$ has holonomy contained in $\Sp{n}$, 
see for instance~\cite{BoGTSM}.

\begin{te}\label{hypcone}
Let $[(H,\Gamma)]$ be a compact locally conformal
hyperk\"ahler manifold. 
Then $H$ is the hyperk\"ahler cone of a compact 3-Sasakian manifold.
\end{te}

\begin{D}
By Theorems~\ref{hypcpt} and \ref{lckcpt} $H$ is the K\"ahler cone of three different compact Sasakian manifold.
Using Theorem~\ref{claimcompact} one sees that they are
isometric, say $W$ the common underlying Riemannian manifold. Then $W$ is
tautologically $3$-Sasakian, because its Riemannian cone $H$ is hyperk\"ahler.
\end{D}

%

\begin{co}
Let $H$ be a hyperk\"ahler manifold. Let $\Gamma$ be a cocompact discrete Lie group 
of hypercomplex homotheties of $H$ acting freely and properly discontinously. 
Then $H$ is the hyperk\"ahler cone of a compact 
3-Sasakian manifold, and $\Gamma$ commutes with the radial flow.
\end{co}

\begin{ex}\label{quathopfsas}
For the quaternionic Hopf manifold given in Example~\ref{quathopf}, 
the 3-Sasakian manifold $W$ given by Theorem~\ref{hypcone} is the standard sphere 
$S^{4n-1}$.
\end{ex}
%
%
%
%
%

\section{Locally conformal hyperk\"ahler reduction}\label{lchkred}

Let $[(H,\Gamma)]$ be a locally conformal hyperk\"ahler manifold.
Denote by $\Homot{\Gamma}{H}$ the Lie group of $\Gamma$-equivariant
hypercomplex homotheties of $H$, and remark that $\Homot{\Gamma}{H}/\Gamma$
coincides with the Lie group $\lchk{M}$ of hypercomplex conformal transformations of 
$M=H/\Gamma$ whenever
$(H,\Gamma)$ is the maximal or the minimal presentation.
Accordingly, a connected Lie subgroup $G$ of $\lchk{M}$ is
{\em twisted  Hamiltonian} if the identity component $\tilde{G}^\circ$ of
$\tilde{G}$ or $G_\text{min}^\circ$ of
$G_\text{min}$ is Hamiltonian in the hyperk\"ahler sense.
We choose the same term Hamiltonian for both the complex and the quaternionic
contexts: the name of Hamilton is of course related both to Hamiltonian
mechanics and to quaternions and it doesn't seem appropriate to look for a
different term for Hamiltonian maps in quaternionic geometry.

If $G$ is a twisted Hamiltonian subgroup of $\lchk{M}$, the quotient of
the hyperk\"ahler momentum map $\tilde{\moment}$ on $\tilde{H}$ is
the {\em locally conformal hyperk\"ahler momentum map}. Nothing different 
if we use the minimal presentation instead.


\begin{te}\label{teonostro}
Let $[(H,\Gamma)]$ be a locally conformal hyperk\"ahler manifold, 
say $M=H/\Gamma$ in the usual sense. Let $G$ be a
twisted Hamiltonian Lie subgroup of $\lchk{M}$, with momentum map $\moment$.
Suppose that the usual topological conditions are satisfied, that is,
$\moment^{-1}(0)$ is non-empty, $0$ is a regular value for $\moment$ and $G$
acts freely and properly on $\moment^{-1}(0)$. Then
\[
[(\tilde{H}\hrid \tilde{G}^\circ,\frac{\Gammamax}{\Gammamax\cap\tilde{G}^\circ})]=[(H_\text{min}\hrid G_\text{min}^\circ,\Gammamin)]
\]
is a locally conformal hyperk\"ahler manifold, where
$\hrid$ denotes the hyperk\"ahler reduction.
\end{te}

\begin{D}
The action of $\tilde{G}^\circ$ commutes with $\Gammamax$ by hypothesis, $\Gammamax$ acts 
properly discontinously and
\[
\frac{\Gammamax}{\Gammamax\cap\tilde{G}^\circ}
\]
acts freely by hypercomplex 
homotheties on the hyperk\"ahler reduction $\tilde{H}\hrid \tilde{G}^\circ$.
The same if we use $(H_\text{min}\hrid G_\text{min}^\circ,\Gammamin)$ instead.
\end{D}

In practice, to reduce a locally conformal hyperk\"ahler manifold $M$ presented as $(H,\Gamma)$ we 
lift $G$ to the hyperk\"ahler universal covering $\tilde{H}$ of $H$, take its identity component
$\tilde{G}^\circ$, and reduce $\tilde{H}$. We then take the hyperk\"ahler reduction 
$\tilde{H}\hrid\tilde{G}^\circ$. The reduction we were looking for is presented by
\[
(\tilde{H}\hrid \tilde{G}^\circ,\frac{\pi_1(M)}{\tilde{G}^\circ\cap\pi_1(M)}).
\]
If we want to get rid of the remaining isometries, because they don't give any contribution
to the locally conformal part of the structure, we switch to the minimal presentation.

\begin{re}
If we start with a subgroup $G$ of $\Homot{\Gamma}{H}$, and we don't know 
it appears as a lifting of a free action on $H/\Gamma$, we must add the hypothesis that 
$\Gamma/(G\cap\Gamma)$ be free and properly discontinous on $H$. 
\end{re}

\begin{de}
Let $[(H,\Gamma)]$ be a locally conformal hyperk\"ahler manifold. Let $G$ be a connected Hamiltonian subgroup 
of $\Homot{\Gamma}{H}$ (this implies that $\Gamma$ is $G$-equivariant) with hyperk\"ahler momentum map \map{\moment}{H}{\liealg{g}^*\otimes\R^3}.
Suppose that $\moment^{-1}(0)$ is non-empty, $0$ is
a regular value for \moment\
and the action of $G$ is free and proper on $\moment^{-1}(0)$, so that the hyperk\"ahler reduction
$H\hrid G$ is defined. Suppose moreover
that the action of $\Gamma$ on $H\hrid G$ be properly discontinous, and that 
the action of $\Gamma/(\Gamma\cap G)$ on $H\hrid G$ be free. Then
\[
[(H\hrid G,\frac{\Gamma}{\Gamma\cap G})]
\]
is the {\em reduction} of the locally conformal hyperk\"ahler manifold $[(H,\Gamma)]$. 
\end{de}

Using Theorem~\ref{hypcone} we relate compact locally conformal
hyperk\"ahler reduction with compact $3$-Sasakian reduction.

\begin{co}\label{trisas}
If $[(\cone{W},\Gamma)]$ is a
compact locally conformal hype\-rk\"a\-hler manifold, and $G$ a twisted
Hamiltonian 
Lie subgroup of $\lchk{M}$ providing a locally conformal hyperk\"ahler reduction, 
then $\tilde{G}^\circ$, $G_\text{min}^\circ$ provides also a
$3$-Sasakian reduction of $\Wmax$, $\Wmin$ respectively.
\end{co}

As in the compact Vaisman case, to any compact locally conformal
hyperk\"ahler manifold $M$ one associates a compact $3$-Sasakian manifold
$W$ as any of the fibers of a canonically defined Riemannian submersion
$M\rightarrow S^1$ and a hypercomplex homothety $\phi$ of the hyperk\"ahler cone of
$W$ such that $M=(\cone{W},\langle \phi \rangle)$ 
\cite[Structure Theorem]{VerVTL}. Accordingly, we have the following analogous of
Corollary~\ref{structure}.

\begin{co}
The Structure Theorem of \cite{VerVTL} and
locally conformal hyperk\"ahler reduction are compatible.
\end{co}

\begin{re}\label{toprem}
Due to Corollary~\ref{trisas}, finding examples of compact locally conformal hyperk\"ahler
reduction accounts in finding discrete Lie groups $\Gamma$ of hypercomplex
homotheties acting freely and properly discontinously on the hyperk\"ahler cone
$\cone{W\hrid G}$ over a reduced $3$-Sa\-sa\-kian manifold $W\hrid G$, and commuting 
with the radial flow of this cone. Here $W$ is a simply connected 3-Sasakian manifold.
According to Proposition~\ref{claimnoncompact}, $\Gamma$ splits as a semidirect product of 
a isometry part $I$ and a ``true'' homothety part $\Z=\langle f\times\rho\rangle$, 
where $f$ is an isometry of $(W\hrid G)/I$ and $\rho$ denotes translations by $\rho\in\R$ 
along the radius of $\cone{(W\hrid G)/I}$.

If we prefer the compact setting, we switch to the minimal presentation. Hence, we look 
for $\Gamma$ acting on $\cone{W\hrid G}$, where $W$ is a compact 3-Sasakian manifold.
In this case, $\Gamma$ splits as $I\rtimes\Z$ for a finite isometry part $I$. 

The isometry type of $\cone{W\hrid G}/\Gamma$ is then related to the isotopy class of $f$ by
\[
\frac{\cone{W\hrid G}}{\Gamma}=\frac{\cone{(W\hrid G)/I}}{\langle f\times\rho\rangle}
=\frac{((W\hrid G)/I)\times[0,\rho]}{(x,0)\sim(f(x),\rho)}.
\]
In particular, $\cone{W\hrid G}/\Gamma$ is a product if and only if $f$ is isotopic to the 
identity. 
\end{re}

\begin{re}
The well known Diamond Diagram pictured in~\cite{BoGTSM}, relating complex, 
positive quaternion-K\"ahler, 3-Sasakian and hyper\-k\"ah\-ler geometries with each other,
has been ported to the compact locally conformal hyperk\"ahler setting in~\cite{OrPLCK}, in the following way.
Let now $(M,[g],I,J,K)$ be 
a compact locally conformal hyperk\"ahler manifold,  where $g$ is the Gauduchon metric
and denote by \lee\ the corresponding length 1 Lee form. Then define 3 foliations 
$\mathcal{B},\mathcal{V},\mathcal{D}$ on $M$, respectively as 
$\langle\lee\rangle,\langle\lee,I\lee\rangle,\langle\lee,I\lee,J\lee,K\lee\rangle$, and suppose the leaves
of these foliations are compact. Then
\begin{itemize}
\item $M/\mathcal{B}$ is a 3-Sasakian orbifold;
\item $M/\mathcal{V}$ is a positive K\"ahler-Einstein orbifold;
\item $M/\mathcal{D}$ is a positive quaternion-K\"ahler orbifold.
\end{itemize}
Moreover, the following diagram holds:
\begin{equation}\label{diamond}
\xymatrix{
 & M\ar[rd]\ar[ld]\ar[dd] & \\
M/\mathcal{V}\ar[rd] & & M/\mathcal{B} \ar[ld]\ar[ll] \\
 & M/\mathcal{D} & \\
}
\end{equation}
Using the argument of Remark~\ref{toprem}, one obtains Diagram~\eqref{diamond} 
from the corresponding DD. For instance, if $M$ is presented as 
$(\cone{W},I\rtimes\Z)$, then $M/\mathcal{B}$ is given by $W/I$. It follows that to any 
reduction of a locally conformal hyperk\"ahler manifold corresponds a reduction of any of 
the spaces in \eqref{diamond}.
\end{re}

\begin{ex}\label{ex1}
Let $M=(\Ha^{n}\setminus 0,\Gamma)$ be a quaternionic Hopf manifold, as described in 
Examples~\ref{quathopf} and \ref{quathopfsas}.
Consider the left action
of $S^1$ on $M$ described by
\begin{equation}\label{action}
t\cdot (h_1,\dots,h_n)=(e^{it}h_1,\dots,e^{it}h_n)\qquad t\in S^1,(h_1,\dots h_n)\in\Ha^{n}\setminus 0.
\end{equation}
This action is $\Gamma$-equivariant and Hamiltonian for the standard metric on $\Ha^{n}\setminus 0$, so
it is by definition a twisted Hamiltonian action on $M$. 
The corresponding 3-Sasakian action (see Corollary~\ref{trisas}) on $S^{4n-1}$ is described by the same Formula
\eqref{action}, where in this case $(h_1,\dots h_n)\in S^{4n-1}$.
The 3-Sasakian
quotient of the sphere is the space (see~\cite{BGMHSS})
\[
\frac{\un{n}}{\un{n-2}\times\un{1}}.
\]
Hence, Remark~\ref{toprem} implies that the locally conformal hyperk\"ahler quotient is
\[
\frac{\un{n}}{\un{n-2}\times\un{1}}\times{\un{1}}.
\]
\end{ex}

\begin{re}
The usual definition of locally conformal hyperk\"ahler structures can be
resumed in the requirement for a hyperhermitian structure with
fundamental forms $\fund_i$ to satisfy $d\fund_i=\lee\wedge\fund_i$ for the 
same closed Lee form \lee. A locally conformal hyperk\"ahler reduction can then be
alternatively defined by means of the Lee form \lee\ and its associated
twisted differential $d^\lee=d-\lee\wedge\cdot$, same way as for the locally conformal
K\"ahler case.
This construction is equivalent with the one given above by means of
presentations.
\end{re}


\begin{re}
In the compact case, we can mimic the proof in \cite{HKLHMS} to provide a direct
construction of the locally conformal hyperk\"ahler quotient. We briefly indicate
the argument. 

Let $(M,[g],I,J,K)$ be a \emph{compact} locally conformal hyperk\"ahler manifold, and denote
by $\fund_i$ the corresponding fundamental forms.
Let $G$ be a twisted
Hamiltonian subgroup of \lchk{M}, and $\moment=(\moment_1,\moment_2,\moment_3)$ 
the corresponding momentum map.
Suppose $0$ is a regular value of $\moment$, and that $G$ acts freely and
properly on $\moment^{-1}(0)$. 
Let $N=\moment_2^{-1}(0)\cap \moment_3^{-1}(0)$. We show that this is a
complex submanifold of $(M,I)$. Indeed, for any $x\in N$, we have
$T_x^\perp N = J\liealg{g}(x)\oplus K\liealg{g}(x)$, where
$\liealg{g}(x)$ is the vector space spanned by the values in $x$ of
the fundamental
fields associated to the elements of $\liealg{g}$. On the other
hand, a vector $X\in T_xM$ belongs to  $T_x\moment_i^{-1}(0)$ if and only if
$\fund_{ix}(X,\cdot)=d^\lee \moment^X_{ix}=0$. Now let $X\in
\Chi(M)$ and $Y$ be any fundamental field. Then, on $N$:
\begin{equation*}
\begin{split}
g(IX,KY)&=-g(JX,Y)=\fund_{2}(X,Y)=0,\\
g(IX,JY)&=g(KX,Y)=\fund_{3}(X,Y)=0
\end{split}
\end{equation*}
Hence $N$ is a complex submanifold of the Vaisman manifold $(M, [g],I)$. 
By construction it is a closed submanifold, thus, by
\cite{TsuHFH} (see also \cite[Proposition 2.1]{VerVTL}), the Lee field is tangent to
$N$ and
$(N,[g_{\mid_N}], I)$ is actually a Vaisman manifold\footnote{The compacity of $M$ is essential to show that $N$ is compact and to apply the cited result.}.
Now $G$ acts on $(N,[g_{\mid_N}], I)$ as a twisted Hamiltonian subgroup of \aut{M}. 
The
associated momentum map is the restriction to $N$ of
$\moment_1$, therefore  $\moment^{-1}(0)=N\cap \moment_1^{-1}(0)$. Then, according
to \cite{GOPLCK}, the
quotient $\moment^{-1}(0)/G$ is a locally conformal K\"ahler manifold with respect to the
projected conformal class and complex structure.
The same argument for $J$ and $K$ ends the proof.
\end{re}


\section{Joyce hypercomplex reduction and HKT reduction}\label{hktred}

Locally conformal hyperk\"ahler manifolds are first of all hypercomplex manifolds, and it
seems natural to ask whether the Joyce hypercomplex reduction given in \cite{JoyHQQ} is
compatible with the locally conformal hyperk\"ahler reduction defined in Section~\ref{lchkred}. 
Recall that in the hypercomplex reduction a momentum
map for the action of a Lie group $G$ on a hypercomplex manifold $(M,I,J,K)$ is
defined as any triple of maps \map{\moment_i}{M}{\liealg{g}^*} satisfying
\begin{equation}\label{joyce}
Id\moment_1=Jd\moment_2=Kd\moment_3
\end{equation}
and $Id\moment_1(X)\neq 0$ for every non-zero fundamental vector field.
%
%
If such
a momentum map exists, then
$(\moment_1^{-1}(0)\cap\moment_2^{-1}(0)\cap\moment_3^{-1}(0))/G$ is a hypercomplex
manifold that we denote by $M\jrid G$.

\begin{te}
Let $G$ be a twisted Hamiltonian Lie subgroup of \lchk{M}. Then a locally conformal
hyperk\"ahler momentum map for the action of $G$ is also a hypercomplex momentum map. Moreover,
\[
M\hrid G=M\jrid G
\]
as hypercomplex manifolds.
\end{te}

\begin{D}
Let $[(H,\Gamma)]$ be the locally conformal
hyperk\"ahler manifold. Denote by $(g,I,J,K)$ the hyperk\"ahler structure on $H$, by $\fund_i$ the fundamental forms and by 
\dismap{\moment_i}{\liealg{g}}{C^\infty(H)}{X}{\moment_i^X\quad i=1,2,3}
the K\"ahler momentum maps. Then by definition
\[
d\moment_i^X=i_X\fund_i
\]
and this implies that
\[
\left.
\begin{split}
Id\moment_1^X&=Ii_X\fund_1\\
Jd\moment_2^X&=Ji_X\fund_2\\
Kd\moment_3^X&=Ki_X\fund_3
\end{split}\right\}=g(X,\cdot).
\]
Hence the locally conformal
hyperk\"ahler momentum map satisfies~\eqref{joyce} and the claim follows.
\end{D}

\begin{re}
A Hermitian manifold equipped with a locally conformal hyperk\"ahler structure might admit Joyce reductions that are not locally conformal hyperk\"ahler reductions.

As an example, let $M=(\Ha^3\setminus 0,\Gamma)$ be the quaternionic Hopf manifold given in Example~\ref{quathopf}, where for simplicity $\Gamma$ is spanned by a real number in $(0,1)$. 
Taking into account that $M$ is diffeomorphic to $S^{11}\times S^1$, we define an action of $S^1\subset\C$ on $M$ by  left multiplication, that is, if $\theta\in S^1$ and $(z_1,\dots,z_7)$ are the complex coordinates describing $S^{11}\times S^1\subset\C^6\times\C$, then $\theta$ acts by

\begin{equation}\label{perandrew}
\theta\cdot (z_1,\dots,z_7)=(e^{2\pi i\theta}z_1,\dots,e^{2\pi i\theta}z_7).
\end{equation}

According to~\cite[Example~6.3]{PePIHS}, this action admits a Joyce momentum map with respect to the standard hypercomplex structure (right multiplication by quaternions).
The hypercomplex quotient is $\su{3}$ which, being simply connected, does not admit any locally conformal hyperk\"ahler structure.
%
 %
%
This means that this Joyce reduction is \emph{not}\/ a locally conformal hyperk\"ahler reduction. In fact this action on $S^{11}\times S^1$ rotates the seventh complex coordinate, that is,  is effective along the $S^1$ factor of the Hopf manifold, hence lifts to an action by non-trivial homotheties on $H^3\setminus0$ which, a fortiori, cannot be a Hamiltonian action.

By passing, we remark the following somehow curious fact (compare also with Example \ref{ex1}).
Denote by $\moment_{\text{Joyce}}$ the lifting to $\Ha^3\setminus 0$ of the Joyce momentum map for the action~\eqref{perandrew} chosen in~\cite{PePIHS}. Consider the following action of $S^1$ on $S^{11}\times S^1$:

\begin{equation}\label{ripetizione}
\theta\cdot (z_1,\dots,z_7)=(e^{2\pi i\theta}z_1,\dots,e^{2\pi i\theta}z_6,z_7).
\end{equation}

This action lifts to the left complex multiplication on $\Ha^3\setminus 0$, which is Hamiltonian with respect to the flat metric. Thus, by definition, the action~\eqref{ripetizione} is twisted Hamiltonian on the quaternionic Hopf manifold and admits a locally conformal hyperk\"ahler momentum map.
Denote by $\moment$ the lifting to $\Ha^3\setminus 0$ of this locally conformal hyperk\"ahler momentum map. Then $\moment_{\text{Joyce}}=\moment$.

Nevertheless, being originated by different actions, this momentum map produces different quotients.
With respect to~\eqref{perandrew}, the quotient is \su{3}, whereas with respect to~\eqref{ripetizione} the quotient is

\[
\frac{\su{3}}{S^1}\times S^1.
\]
%
%
\end{re}

%
%
%
%
%
%
%
%
%
%

We now discuss the relation between locally conformal hyperk\"ahler reduction and HKT reduction. Recall first that on a hyperhermitian manifold $(M,g)$, a connection $D$ is called hyperk\"ahler with torsion (HKT), if it is hyperhermitian (namely, $DI=DJ=DK=0$, $Dg=0$ and its torsion is totally skew-symmetric (that is, $g(\Tor^D(X,Y),Z)$ is a $3$-form).

The following result shows that compact locally conformal hyper\-k\"ah\-ler geometry (in which the
Lee form can be always assumed parallel and unitary) is related to HKT
geometry.

\begin{te}\cite{OPSPFH}
Let $(M, [g], I,J,K)$ be a locally conformal hyperk\"ahler manifold with unitary parallel Lee form
$\lee$ associated to the metric $g$. Then the metric
\begin{equation}\label{hkt}
\hat g=g-\frac 12\{\lee\otimes\lee+
I\lee\otimes I\lee+J\lee \otimes J\lee+
K\lee\otimes K\lee\}
\end{equation}
is HKT.
\end{te}

As a  reduction scheme was recently provided for HKT structures (see~\cite{GPPRHS}),
it is natural to discuss the interplay between the two reductions.

Recall that, based on Joyce's hypercomplex reduction (see~\cite{JoyHQQ}), the HKT reduction has the peculiarity that
a ``good'' action of a group of hypercomplex automorphisms does not automatically
produce a momentum map. What is proved in \cite{GPPRHS} is that if a group acts by
hypercomplex isometries with respect to the HKT metric and {\em if Joyce's
hypercomplex momentum map exists}, then the hypercomplex quotient inherits a
natural HKT structure.

\begin{te}\label{comphkt}
Let $M$ be a compact locally conformal hyperk\"ahler manifold, $G$ a compact
twisted Hamiltonian Lie group and \moment\ the corresponding momentum map, with
$0$ as a regular value.
Suppose that $G$ acts freely on $\moment^{-1}(0)$, so that
the locally conformal hyperk\"ahler reduction of Theorem~\ref{teonostro} is defined. Then 
\moment\ is a hypercomplex momentum map, and $G$ acts by isometries for the 
associated HKT metric. Moreover,
the HKT structure of the quotient is induced by the locally conformal hyperk\"ahler 
reduced structure by the relation \eqref{hkt}.  
\end{te}

%
\begin{D}
Let $g$ be the Gauduchon metric of $M$. Then $G$ acts by isometries for $g$, thus preserving
the corresponding Lee form, and since $\hat g$ is defined using only $G$-invariant tensors, then 
$G$ acts also by $\hat g$-isometries. In this situation \moment\  
satisfies Joyce's conditions, hence
\moment\  is also a HKT momentum map and the HKT quotient exists. Moreover, the
reduced HKT structure lives on the reduced locally conformal hyperk\"ahler manifold. 
As the projection
$\moment^{-1}(0)\rightarrow \moment^{-1}(0)/G$ is a Riemannian submersion with respect to
both the locally conformal hyperk\"ahler and HKT metric, and since the Lee form is projectable, 
we see that
the reduced HKT metric is the one induced by the reduced locally conformal hyperk\"ahler
metric.
\end{D}

%
Compact HKT manifolds don't have a potential, but if the HKT structure is induced
by a locally conformal hyperk\"ahler structure, then the Lee form is
a \emph{potential 1-form} as defined in~\cite{OPSPFH}. Thus we have the following
Corollary.

\begin{co}
All manifolds obtained by the HKT reduction associated to a locally conformal 
hyperk\"ahler reduction admit a potential 1-form.
\end{co}

\begin{re} 
It is not known if generic HKT reduction preserves potential 1-forms.
\end{re}

\begin{re} 
Not all HKT structures are induced by a locally conformal hyperk\"ahler structure,
see~\cite{VerHMT}, so one may ask if the HKT reduction can be induced by a 
locally conformal hyperk\"ahler structure in such cases.
\end{re}


\begin{thebibliography}{HKLR87}

\bibitem[Bel00]{BelMSN}
F.~Belgun.
\newblock {On the metric structure of the non-K\"ahler complex surfaces}.
\newblock {\em Math. Ann.}, 317:1--40, 2000.

\bibitem[BG99]{BoGTSM}
C.~Boyer and K.~Galicki.
\newblock 3-{S}asakian manifolds.
\newblock In {\em Surveys in differential geometry: essays on Einstein
  manifolds}, pages 123--184. Int. Press, Boston, MA, 1999.

\bibitem[BGM98]{BGMHSS}
C.~P. Boyer, K.~Galicki, and B.~M. Mann.
\newblock {Hypercomplex structures from 3-Sasakian structures}.
\newblock {\em J. Reine Angew. Math.}, 501:115--141, 1998.

\bibitem[CP99]{CaPEWG}
D.~M.~J. Calderbank and H.~Pedersen.
\newblock Einstein-{W}eyl geometry.
\newblock In {\em Surveys in differential geometry: essays on Einstein
  manifolds}, pages 387--423. Int. Press, Boston, MA, 1999.

\bibitem[CS04]{CaSAHS}
F.~M. Cabrera and A.~Swann.
\newblock {Almost Hermitian structures and quaternionic geometries}.
\newblock {\em Differ. Geom. Appl.}, 21(2):199--214, 2004.

\bibitem[Gau95]{GauSWE}
P.~Gauduchon.
\newblock Structures de {W}eyl-{E}instein, espaces de twisteurs et vari\'et\'es
  de type {$S\sp 1\times S\sp 3$}.
\newblock {\em J. Reine Angew. Math.}, 469:1--50, 1995.

\bibitem[GO98]{GaOLCK}
P.~Gauduchon and L.~Ornea.
\newblock Locally conformally {K}\"ahler metrics on {H}opf surfaces.
\newblock {\em Ann. Inst. Fourier}, 48:1107--1127, 1998.

\bibitem[GOP05]{GOPLCK}
R.~Gini, L.~Ornea, and M.~Parton.
\newblock {Locally conformal K\"ahler reduction}.
\newblock {\em J. Reine Angew. Math.}, 581, April 2005.

\bibitem[GPP02]{GPPRHS}
G.~Grantcharov, G.~Papadopoulos, and Y.~S. Poon.
\newblock Reduction of {HKT}-structures.
\newblock {\em J. Math. Phys.}, 43(7):3766--3782, 2002.

\bibitem[HKLR87]{HKLHMS}
N.~J. Hitchin, A.~Karlhede, U.~Lindstr{\"o}m, and M.~Ro{\v{c}}ek.
\newblock Hyper-{K}\"ahler metrics and supersymmetry.
\newblock {\em Comm. Math. Phys.}, 108(4):535--589, 1987.

\bibitem[Ino74]{InoSCV}
M.~Inoue.
\newblock {On surfaces of class VII$_0$.}
\newblock {\em Invent. Math.}, 24:269--310, 1974.

\bibitem[IPP05]{IPPLCP}
S.~Ivanov, M.~Parton, and P.~Piccinni.
\newblock {Locally conformal parallel $\text{\upshape \rmfamily G}_2$,
  $\text{\upshape \rmfamily Spin}(7)$ and $\text{\upshape \rmfamily Spin}(9)$
  structures}.
\newblock In preparation, 2005.

\bibitem[Joy91]{JoyHQQ}
D.~Joyce.
\newblock {The hypercomplex quotient and the quaternionic quotient}.
\newblock {\em Math. Ann.}, 290(2):323--340, 1991.

\bibitem[KO05]{KaOGFC}
Y.~Kamishima and L.~Ornea.
\newblock {Geometric flow on compact locally conformally {K}\"ahler manifolds}.
\newblock {\em Tohoku Math. J., II. Ser.}, 57(2), June 2005.

\bibitem[Kod64]{KodSC1}
K.~Kodaira.
\newblock On the structure of compact complex analytic surfaces, {I}.
\newblock {\em American J. Math.}, 86:751--798, 1964.

\bibitem[Kod66]{KodSC2}
K.~Kodaira.
\newblock On the structure of compact complex analytic surfaces, {II}.
\newblock {\em American J. Math.}, 88:682--721, 1966.

\bibitem[Kod68a]{KodSC3}
K.~Kodaira.
\newblock On the structure of compact complex analytic surfaces, {III}.
\newblock {\em American J. Math.}, 90:55--83, 1968.

\bibitem[Kod68b]{KodSC4}
K.~Kodaira.
\newblock On the structure of compact complex analytic surfaces, {IV}.
\newblock {\em American J. Math.}, 90:1048--1066, 1968.

\bibitem[OP97]{OrPLCK}
L.~Ornea and P.~Piccinni.
\newblock {Locally conformal K{\"a}hler structures in quaternionic geometry}.
\newblock {\em Trans. Am. Math. Soc.}, 349(2):641--655, 1997.

\bibitem[OPS03]{OPSPFH}
L.~Ornea, Y.~S. Poon, and A.~Swann.
\newblock {Potential 1-forms for hyper-K\"ahler structures with torsion}.
\newblock {\em Classical Quantum Gravity}, 20(9):1845--1856, 2003.

\bibitem[OT05]{OeTNKC}
K.~Oeljeklaus and M.~Toma.
\newblock {Non-K\"ahler compact complex manifolds associated to number fields}.
\newblock {\em Ann. Inst. Fourier}, 55(1):1291--1300, 2005.

\bibitem[OV03]{OrVSTC}
L.~Ornea and M.~Verbitsky.
\newblock Structure theorem for compact {V}aisman manifolds.
\newblock {\em Math. Res. Lett.}, 10(5-6):799--805, 2003.

\bibitem[PP99]{PePIHS}
H.~Pedersen and Y.-S. Poon.
\newblock Inhomogeneous hypercomplex structures on homogeneous manifolds.
\newblock {\em J. Reine Angew. Math.}, 516:159--181, 1999.

\bibitem[PPS93]{PPSEWE}
H.~Pedersen, Y.~S. Poon, and A.~Swann.
\newblock The {E}instein-{W}eyl equations in complex and quaternionic geometry.
\newblock {\em {D}ifferential {G}eometry and its {A}pplications},
  3(4):309--322, 1993.

\bibitem[Tsu94]{TsuHFH}
K.~Tsukada.
\newblock {Holomorphic forms and holomorphic vector fields on compact
  generalized Hopf manifolds}.
\newblock {\em Compos. Math.}, 93(1):1--22, 1994.

\bibitem[Vai79]{VaiLCK}
I.~Vaisman.
\newblock Locally conformal {K}\"ahler manifolds with parallel {L}ee form.
\newblock {\em Rend. Mat. Roma}, 12:263--284, 1979.

\bibitem[Vai82]{VaiGHM}
I.~Vaisman.
\newblock Generalized {H}opf manifolds.
\newblock {\em Geometriae Dedicata}, 13:231--255, 1982.

\bibitem[Ver03]{VerHMT}
M.~Verbitsky.
\newblock {Hyperk\"ahler manifolds with torsion obtained from hyperholomorphic
  bundles}.
\newblock {\em Math. Res. Lett.}, 10(4):501--513, 2003.

\bibitem[Ver04]{VerVTL}
M.~Verbitsky.
\newblock {Theorems on the vanishing of cohomology for locally conformally
  hyper-K\"ahler manifolds. (Russian)}.
\newblock {\em Tr. Mat. Inst. Steklova, Algebr. Geom. Metody, Svyazi i
  Prilozh}, (246):64--91, 2004.
\newblock arXiv:math.DG/0302219.

\end{thebibliography}
\end{document}